\documentclass[11pt]{article}

\usepackage{amsmath,makeidx,amssymb,amsthm,amscd,a4}
\usepackage{graphicx,epsfig,latexsym,index}
\newtheorem{theor}{Theorem}
\newtheorem{lemma}[theor]{Lemma}

\newcommand{\me} {\tiny \mbox{\raisebox{.6ex}{1 \hspace{-2.2mm} 
\raisebox{-.4ex}{/ \hspace{-2.2mm} \raisebox{-.4ex}{2} } } } \hspace{-2.6mm} }
\newcommand{\R}{\mathbb{R}}
\newcommand{\eps}{\epsilon}
\newcommand{\vphi}{\varphi}
\begin{document}


\title{Applications of the Backus--Gilbert method to linear and some 
non--linear equations}

\author{A. Leit\~ao \\ \small
Dep. of Mathematics, Federal Univ.\,of St.\,Catarina, 88040-900 Florian\'opolis, Brazil}

\maketitle

\begin{abstract}
        We investigate the use of a functional analytical version of the 
Backus--Gilbert Method as a reconstruction strategy to get specific 
information about the solution of linear and slightly non-linear systems 
with Frech\'et derivable operators. Some {\em a priori} error estimates are 
shown and tested for two classes of problems: a non-linear {\em moment 
problem} and a linear elliptic {\em Cauchy problem}. For this second class 
of problems a special version of the Green-formula is developed (see 
Theorem~\ref{korrektur_formel}), in order to analyze the involved adjoint 
equations.
\end{abstract}

\section{Introduction}

\subsection{Main results}

        The functional analytical approach of the Backus--Gilbert method in 
\S 1.3 was already used by other authors (see [Ch], [Ki] or [LM1,2]). 
In this paper we use the differentiability of the involved non-linear 
operator in order to develop the error estimative (\ref{err_bg_nicht_lin}) 
and (\ref{err_nicht_lin_reg}) for this reconstruction schema. If the operator 
is linear, we obtain the estimative (\ref{err_bg_lin}).

        In order to test this reconstruction strategy, we choose the same 
non-linear operator in \S 3.1 as Louis does in [Lo]. The numerical tests in 
\S 4.1 show that one can get good results even for noisy data.

        For the second test in \S 3.2, we choose a linear operator, which is 
highly ill-posed. The results are again satisfactory provided one uses 
apropriated {\em sentinels} to define the reconstruction strategy (see 
[Ch] or [Le]).

        The results presented in this paper constitute part of the author's 
PhD research and they can be found with a bit more detail in [Le].

\subsection{Historical overview}

        This reconstruction method was first proposed in 1967 by G.Backus and 
J.Gilbert [BG1,2,3]. They were interested in the pointwise reconstruction 
of a function $f\in X = L^2(\Omega)$, were $\Omega\subset \R^n$ is supposed 
to be open and bounded. The motivation of their problems was geophysical 
and the mathematical problem involved in the model is known in the 
literature as the {\em moment problem}. It can be formulated as follows: 
find a function $f \in X$ such that
\begin{equation}
        \int_{\Omega}{\, K_i(x)\, f(x)\ \d x} \ = \ g_i \, ,\ \ i=1,\dots,N,
\label{moment_probl} \end{equation}
were the kernels $K_i$ are known real functions, which are well defined 
at $\Omega$ and the right hand side $g = \{g_i\}_{i=1}^{N} \in Y = \R^N$ 
correspond to the measured data of the physical problem. In order to 
determine the value of the solution $f$ at some point $x_0 \in \Omega$, 
they suggested a linear reconstruction schema, which is defined by a 
functional of the right hand side of the linear system (\ref{moment_probl}). 
One defines the linear functional $R_N \in Y'$ by 
\begin{equation} 
    R_N(g) \ := \ < \varphi, g >_Y \
           \  = \ \int_{\Omega}{\underbrace{ \left( \sum_{i=1}^N
              {\varphi_i K_i(x)} \right)}_{\phi_N(x)}f(x) \ \d x}
           \  = \ < \phi_N, f >_X,
\end{equation}
were $\varphi \in Y'$ and $\phi_N \in X$. It is easy to observe that 
$f_N(x_0) := R_N(g)$ will be a good approximation for $f(x_0)$ if the 
condition $\phi_N(\cdot) \simeq \delta(x_0-\cdot)$ is satisfied. The 
Backus--Gilbert idea is to force this condition by defining the quadratic 
functional
\begin{equation}
             J(\phi) \ := \ \int_\Omega{\, |x_0-x|^2\, \phi^2(x) \ dx}
\label{BG_funktional} \end{equation}
on $X$ and choosing $\phi_N$ such that
\begin{equation}
          J(\phi_N) \ = \ \min_{\phi \, \in \, {\rm Span}\{K_i\}}{\, J(\phi)}.
\label{BG_minimierung} \end{equation}
The linear constraint
\begin{eqnarray}
                  \int_{\Omega}{\phi(x) \ dx} \ = \ 1  \nonumber
\end{eqnarray}
is imposed in order to avoid the trivial solution in (\ref{BG_minimierung}). 
Once one has evaluated the function $\phi_N(x) = \sum{\varphi_i K_i(x)}$, 
the approximation $f_N(x_0)$ is determined by the inner product
\begin{equation}
            f_N(x_0) \ = \ < \phi_N, f >_X \ = \ < \varphi, g >_Y,
\end{equation}
were $\varphi = \{\varphi_i\}_{i=1}^N$. One great advantage of using the 
Backus--Gilbert method which can be recognized in (\ref{BG_minimierung}) 
is that the evaluation of the reconstruction operator $R_N(\cdot)$ does 
not depend on the system data. For different sets of data $g$ is possible 
to reconstruct the value of the respective $f(x_0)$ only by evaluating an 
inner product in $Y$.

\subsection{Functional analytical formulation}

        Let $V\hookrightarrow X\hookrightarrow V'$ be a Hilbert triple, 
$Y$ a Hilbert space, $y \in Y$ and $A: X \rightarrow Y$ a bounded linear 
operator. We analyze the problem of finding the value $< \mu, x^*>$ for 
$\mu \in V'$, were $x^*$ is the generalized solution obtained by the 
Moore--Penrose inverse of
\begin{equation}                A \, x \ = \ y.
\label{allg_moment_probl} \end{equation}
It is obvious that the expression $<\mu, x^*>$ does not need to be well 
defined, if we do not make any further regularity assumptions about 
$x^*$. Depending on the physical situation involved, it is possible to 
guarantee that the expression $<\mu, x^*>$ is well defined for some $\mu$'s 
or even that $x^* \in V$. As we suppose $y$ is obtained by measurements, it 
is to be expected that only a $y_\eps$ with $||y - y_\eps||_Y \leq \eps$ is 
available, with $\eps > 0$ small.

        We use the Backus--Gilbert strategy and try to reconstruct the value 
$f:=$ $<\mu, x^*>_X$ using a linear functional evaluated in $y_\eps$. 
For $\vphi \in Y'$ we define $f_{\eps,\vphi} := <\vphi, y_\eps>_Y$ 
and estimate the error $|f - f_{\eps,\vphi}|$ by
\begin{eqnarray}
|f-f_{\epsilon,\varphi}| & = &
                      |<\mu, x^*>_X - <\vphi, y_\eps>_Y| \nonumber \\
  & \leq & |<\vphi, y - y_\eps>_Y| \ + \
                            |<\mu, x^*>_X - <\vphi, A\, x^*>_Y| \nonumber \\
  & \leq & \eps\, ||\vphi||_{Y'}\ +\ |<\mu - A^*\, \vphi, x^*>_X|,
                                                         \label{err_f_f.ep.fi}
\end{eqnarray}
were $A^*:Y'\rightarrow X'$ is the adjoint operator of $A$. If we succeed 
in finding a solution $\vphi \in Y'$ for the equation $A^* \vphi = \mu$ 
we can write
\begin{equation}
f \ = \ <A^*\, \vphi, x^*>_X \ = \ <\vphi, y>_Y \ \simeq \
        <\vphi, y_\eps>_Y \ = \ f_{\eps, \vphi},
\label{aprox_f_f.ep.fi} \end{equation}
and the error $|f - f_{\eps,\vphi}|$ behaves like $O(\eps)$. Another 
consequence is that the approximation $f_{\eps,\vphi}$ is exact if 
there are no errors in the measurements ($y_\eps=y$).

        In the special case of $X$ and $Y$ being spaces of functions defined 
over a region $\Omega$, the Backus--Gilbert strategy suggests a 
pointwise reconstruction of $x^*$. In order to reconstruct the value of 
$f(\cdot)$ at the point $t \in \Omega$ we should take $\mu(\cdot) = 
\delta(t-\cdot)$ in (\ref{aprox_f_f.ep.fi}) and solve the adjoint equation 
$A^* \vphi = \delta$.%
\footnote{The Hilbert space $V$ must be chosen, so that $\delta$ belongs 
to $V'$.}

        We may have difficulties if $\delta \not\in \overline{Rg(A^*)}$. 
In this case we can use the projection of $\delta$ over $Ker(A^*)^\perp$ 
instead of $\delta$ itself. This is equivalent to minimizing the error 
$||A^*\, \varphi - \delta||^2_{V'}$ or to find a solution $\varphi \in Y'$ 
of the normal equation
\begin{eqnarray}
                (A\, A^*) \, \varphi \ = \ A\, \delta.  \nonumber
\end{eqnarray}

        Louis and Maa{\ss} propose a similar approach in [LM2] and use the 
projection of $\delta$ over special Sobolev spaces of negative index. In 
[LM1] (see also [Lo]) the equation $(A \, A^*) \, \varphi = A \, e_h$ is 
considered, were $e_h$ is a {\em mollifier}, i.e., a smooth approximation 
for the Dirac-distribution $\delta$.

        Another alternative for the case $\delta \not\in Rg(A^*)$ was 
proposed by Chavent in [Ch]. He tried to regularize the normal equations 
using the Tikhonov-strategy: $\varphi$ is chosen as the minimum over $Y'$ 
of the functional $(\, ||A^* \varphi - \delta||_{V'}^2 + \alpha 
||\varphi||_{Y'}^2\, )$, were $\alpha > 0$ is a small regularization 
parameter.

\section{Analysis of the method}

        We are interested in applying the Backus--Gilbert strategy for 
operators of the form $A = A_0 + \gamma A_1$, were $A_0 \in {\cal L}(X,Y)$, 
\ $A_1: X \mapsto Y$ is continuously differentiable%
\footnote{The Fr\'echet derivative of $A_1$ will be denoted by $d A_1$.}
in $X$ and $\gamma > 0$ is a small number. Let $\mu \in V'$ and $y_\eps 
\in Y$ as before.%
\footnote{For convenience we will identify the spaces $X$ with $X'$ and $Y$ 
with $Y'$.}

\begin{lemma}

        If $x^0 \in V$ is an approximation to a solution $x^*$ of 
(\ref{allg_moment_probl}), the expression $f_{\eps,\varphi} = <\vphi, 
y_\eps>_Y$ gives an approximation for $f := <\mu, x^* >_X$ and the error 
$|f - f_{\eps,\vphi}|$ is estimated by
\begin{eqnarray} 
|f-f_{\eps,\vphi}| & \leq & |<\vphi,y-y_\eps>_Y| \ +\
        |<\vphi,A\, x^* - A\, x^0 - dA(x^0)(x^*-x^0)>_Y|  \nonumber \\[0.1cm]
  & & + \ |<\vphi, A\, x^0 - dA(x^0)x^0>_Y|\ +\
                    |<dA(x^0)^*\vphi - \mu,x^*>_X|. \label{err_bg_nicht_lin}
\end{eqnarray}
\end{lemma}
{\it Proof:}  Estimate (\ref{err_bg_nicht_lin}) follows promptly from the 
following equality
\begin{eqnarray*}
|f-f_{\eps,\vphi}| &=& |<\mu, x^*>_X - <\vphi, y_\eps>_Y \pm \\[0.1cm]
                   & & \pm <\vphi, y> \pm <\vphi, A\, x^0 - dA(x^0)(x^*-x^0)>|
\end{eqnarray*}
\hfill 

        Before analyzing the right hand side of (\ref{err_bg_nicht_lin}), 
let us discretizate the spaces involved. We define the finite dimensional 
space $Y_h =$ Span$\{y_j\}_{j=1}^N$ by
\begin{equation}
   Y_h\ \subset\ \{ \varphi\in Y\ /\ <\varphi,A\, x^0-dA(x^0)x^0>_Y = 0 \}.
\label{Yh_wahl} \end{equation}
Further we let $P_h: Y \mapsto Y_h$ be the orthogonal projector over $Y_h$ 
and choose the finite dimensional space $X_h =$ Span$\{x_j\}_{j=1}^N\ \subset 
D(A) \cap V$ such that the property
\begin{equation}
       \det{(<dA(v)^* P_h^* y_i, x_j>)}_{1\leq i,j\leq N} \ \not= \ 0.
\label{det_Xh_Yh} \end{equation}
is satisfied.

\begin{theor}

        Define $y_h:= P_h \, y$, \ $y_{\eps,h}:= P_h \, y_\eps$  and \ 
$f_{\eps,\vphi,h}:= \, <\vphi, y_{\eps,h}>_Y$. For every $\vphi \in Y$ the 
following estimate holds:
\begin{eqnarray}
|f-f_{\eps,\vphi,h}| & \leq & \eps \ ||P_h|| \ ||\varphi||_Y  +
                        \gamma \ ||\vphi||_Y \ ||P_h|| \ O (||x^* - x^0||_X^2)
                                                           \nonumber \\[0.1cm]
   & & + |< dA^*(x^0) P_h^* \vphi - \mu , x^* >_X|. \label{err_nicht_lin_absc}
\end{eqnarray}
\end{theor}
{\it Proof:} By an argument analogous to that used in (\ref{err_bg_nicht_lin}) 
we obtain that for each $\vphi \in Y$
\begin{eqnarray} 
|f-f_{\eps,\vphi,h}| & \leq &
  ||\vphi||_Y \, ||y_h - y_{\eps,h}||_Y \nonumber \\
& & \ + \ \gamma\, ||\vphi||_Y \, ||P_h\, A_1\, x^* - P_h\, A_1\, x^0
                                 - P_h\, dA_1(x^0)(x^* - x^0)||_Y \nonumber \\
& & \ + \ |< \vphi , P_h\, A\, x^0 - P_h\, dA(x^0) x^0 >_Y| \nonumber \\
& &\ +\ |< dA^*(x^0) P_h^* \vphi - \mu , x^* >_X|. \label{err_nicht_lin_reg}
\end{eqnarray}
The first term in (\ref{err_nicht_lin_reg}) can be estimated by
$$  ||\vphi||_Y \, ||y_h - y_{\eps,h}||_Y \ \leq \
                                           \eps \ ||P_h|| \ ||\varphi||_Y.  $$
For the second term we have
$$  \gamma\, ||\vphi||_Y \, ||P_h\, A_1\, x^* - P_h\, A_1\, x^0
                              - P_h\, dA_1(x^0)(x^* - x^0)||_Y \ \leq \
                   \gamma \ ||\vphi||_Y \ ||P_h|| \ O (||x^* - x^0||_X^2).  $$
\noindent  The third term in (\ref{err_nicht_lin_reg}) disappears because of 
our choice of $Y_h$. Putting these inequalities together we obtain 
(\ref{err_nicht_lin_absc}). \hfill 

        The last term in (\ref{err_nicht_lin_absc}) gives us a rule for 
choosing our $\varphi \in Y$. This is actually
\begin{equation}
  < dA^*(x^0) P_h^* \varphi, \ x_j >_X \ = \ <\mu, x_j>_X,\ \ j=1 \dots N.
\label{allg_adjun_gl} \end{equation}
That means we can evaluate the coefficients of $P_h \, \varphi$ in $Y_h$ 
by solving the N-dimensional linear system (\ref{allg_adjun_gl}). Solving 
this system is a well defined problem, as can be seen from the determinant 
condition (\ref{det_Xh_Yh}).

        Next we interpret the system (\ref{allg_adjun_gl}) in a different 
way. Let us assume that the space $X_h$ can be written as $B^* \, Y_h$, 
where $B$ is a linear bounded operator $B: V' \mapsto Y$ with $B^*: Y 
\mapsto V$. We are then able to write (\ref{allg_adjun_gl}) as
\begin{eqnarray*}
  < dA^*(x^0) P_h^* \, \vphi - \mu , \ B^* w >_X \ = \ 0,\
                                                         \ \forall w\in Y_h,
\end{eqnarray*}
i.e.,
\begin{equation}
  < B G_h \vphi , \ w >_Y \ = \ < B \mu , \ w >_Y,\ \ \forall w\in Y_h,
\label{bg_projektion} \end{equation}
\noindent  where $G_h = dA^*(x^0) P_h^*$. If $\mu = \delta$ and we are in the 
special case $\mu \in$ Ker$\,B$, it follows from (\ref{bg_projektion}) that
\begin{equation}
           < B G_h \vphi , \ w >_Y \ = \ 0,\ \ \forall w\in Y_h.
\label{bg_proj_spez} \end{equation}
\noindent  Further if it is possible to decompose the product $B G_h$ as a 
square ${\cal B}^2$ of a symmetric matrix ${\cal B}$, it follows from 
(\ref{bg_proj_spez}) that $\| {\cal B} \vphi \|^2_Y = 0$. Instead of solving 
system (\ref{bg_proj_spez}), we can consider the minimization problem:
$$  \left\{ \renewcommand{\arraystretch}{1.5} \begin{array}{l}
        ||{\cal B} \vphi_h||^2 \ = \ \min_{\, \vphi \in Y_h}
                                       {\, ||{\cal B} \vphi||^2} \\
        {\rm under\ the\ linear\ constraint}\ \ < dA^*(x^0) \vphi_h , 1 >_X \
                                                                        = \ 1
    \end{array} \right. $$
The extra linear constraint is motivated by the original Backus--Gilbert 
formulation in \S 1.2 and introduced in order to avoid the trivial solution 
in the minimization problem. The constrained minimization problem above can 
be interpreted as an extended Backus--Gilbert method.

        We proceed to develop an error estimate for the linear case when 
noisy data is considered.

\begin{theor}

        Let $A$ be a linear operator. Take $B^* = A^*$ and $X_h = A^* Y_h$. 
If we choose $\vphi_h \in Y_h$ to be the solution of (\ref{allg_adjun_gl}), 
i.e., $<A^* \vphi_h - \mu, w> = 0, \, \forall w \in X_h$, we obtain the 
error estimate
\begin{equation} 
|f - f_{\eps,\vphi,h}| \ \leq \ {\rm dist}(\mu , A^* Y_h)_{V'} \
                (1 + ||{\cal P}_h||) \ {\rm dist}(x^* , X_h)_V \ + \ O (\eps).
\label{err_bg_lin} \end{equation}
\end{theor}
{\it Proof:}  Using (\ref{err_f_f.ep.fi}), for each $\vphi_h = P_h \vphi 
\in Y_h$ we obtain the equality
\begin{eqnarray*}
|f - f_{\eps,\vphi,h}| \ = \ O (\eps) \ + \ |< A^* \vphi_h - \mu , x^* >_X|.
\end{eqnarray*}
Define $x_h := {\cal P}_h x^*$, were ${\cal P}_h: X \mapsto X_h$ is 
the orthogonal projector over $X_h$. Now choosing $\vphi_h \in Y_h$ as 
the solution of (\ref{allg_adjun_gl}), for every $\psi \in Y_h$ we have
\begin{equation} 
<\! A^* \vphi_h - \mu , x^* \!>_X \ = \
                     <\!A^* \vphi_h - A^* \psi, x^* - x_h\!>_X \, + \,
                     <\!A^* \psi - \mu, x^* - x_h\!>_X. \label{manipul_fehler}
\end{equation}
The first term on the right hand side of (\ref{manipul_fehler}) disappears 
by the definition of $x_h$. For the second term we have
\begin{eqnarray*}
  |<\!A^* \varphi_h - \mu ,\ x^* - x_h\!>_X| \ \leq \ 
                                  ||A^* \psi_h - \mu||_{V'} \ ||x^* - x_h||_V.
\end{eqnarray*}
In order to estimate $||x^* - x_h||_V$ we define $\tilde{x} \in 
X_h$ as the solution of the minimization problem
\begin{eqnarray*}
       ||x^* - \tilde{x}||_V^2 \ = \ \min_{x \in X_h}{\, ||x^* - x||_V^2}.
\end{eqnarray*}
From this definition follows
\begin{eqnarray*}
||x^* - x_h||_V &\leq & ||x^* - \tilde{x}||_V\ +\
                                            ||\tilde{x} - {\cal P}_hx^*||_V \\
& \leq & {\rm dist}(X_h, x^*)_V \ + \ ||{\cal P}_h(\tilde{x} - x^*)||_V \\
& \leq & (1 + ||{\cal P}_h||) \ {\rm dist}(X_h, x^*)_V,
\end{eqnarray*}
and the teorem is prooved. \hfill 

        It is easy to conclude from (\ref{err_bg_lin}) that the error in the 
approximation $f_{\eps,\vphi,h}$ will converge to zero with $h$ and $\eps$ 
only when we have $\mu \in \overline{Rg\, A^*}$.

\section{Applications}

\subsection{A non-linear moment problem}

        We start this discussion with a special class of non-linear 
moment problems. Quadratic moment problems were also analyzed by Louis in 
[Lo]. Let $X = Y = L^2(0,1)$ and $A:X \mapsto Y$ the operator defined by
\begin{equation} 
    (Ax)(t) \ = \ \int_0^t {x^0(t-s)\ x(s)\ ds} \ + \
                          \nu \int_0^t {x(t-s)\ x(s)\ ds}, \ \ t\in [0,1],
\label{nicht-lin_operator} \end{equation}
\noindent  were the kernel $x^0$ of the linear component of $A$ is a 
$L^2(0,1)$-function and $\nu > 0$ is a small parameter, that controls the 
non-linear component of $A$. Just like in \S 2 we will analyze the 
system \ $A x \, = \, (A_0 + \nu A_1)\, x \, = \, y$.

        The right hand side of this system consists of measured data, so we 
assume we know only a finite number of $y_i = y(t_i)$, $t_i \in (0,1)$. If we 
define the projection operator $P_h: Y \mapsto Y_h = \R^N$, it is possible to 
define a discrete version of $A$ in (\ref{nicht-lin_operator}) by setting
\begin{eqnarray*} 
  (P_h A)(x)\ :=\ \left[ Ax (t_i) \right]^t \ = \ 
                    \left[ \ \int_0^{t_i}{x^0(t_i-s)\ x(s)\ ds} \ + \
                        \nu \int_0^{t_i}{x(t_i-s)\ x(s)\ \d s} \ \right]^t.
\end{eqnarray*}
\noindent  If we further assume that our measurements are inexact, then we 
actually have a $y_{h,\eps} \in Y_h$ with $||P_h y - y_{h,\eps}|| \leq 
\eps$, were $\eps > 0$ is small. We will be interested in finding the 
solution $x^*$ of the discrete non-linear system
\begin{eqnarray*}
                           (P_h A)\, x^* \ = \ y_{h,\eps}.
\end{eqnarray*}
        We saw in \S 2 that an approximation for the solution $x^*$ 
is needed. For this propose we will choose the kernel $x^0$ in 
(\ref{nicht-lin_operator}). We also need the operators $P_h d A$ and its 
adjoint $(P_h d A)^*:Y_h \mapsto X$. One can easily see that for $f \in 
L^2(0,1)$ and $w \in \R^N$ the equalities
\begin{eqnarray*} 
  (P_h d A(x))(f) \ = \ \left[ \ \int_0^{t_i}{x^0(t_i-s)\ f(s)\ \d s}\ +\
         2\nu \int_0^{t_i}{x(t_i-s)\ f(s)\ \d s} \ \right]^t_{1\leq i\leq N}
\end{eqnarray*}
\noindent and
\begin{eqnarray*}
 (P_h d A(x))^*(w) \ = \ \sum_{i=1}^N {\, w_i \ [x^0(t_i-s) + 
                                      2\nu x(t_i-s)] \ \chi_{[0,t_i]}(s)}.
\end{eqnarray*}
\noindent  are valid.

        Now we have to choose the space $X_h =$ Span$\{x_i\}_{i=1}^N$. This 
choice must reflect the expected regularity of the solution $x^*$ and 
should be such that the system in (\ref{allg_adjun_gl}) has nice properties. 
We choose a cubic B--spline basis for $X_h$ for the numerical experiments. 
Given $\mu \in X'$ we will have to solve the system
\begin{equation}
\left[ < (P_h d A(x^0))^* e_i , x_j >_X \right]_{i,j=1}^N \ [ \vphi_j ]^t
\ = \ \left[ < \mu , x_j >_X \right]^t,
\label{nicht_lin_moment_probl_syst} \end{equation}
\noindent  where the matrix of (\ref{nicht_lin_moment_probl_syst}) will have 
almost upper triangular form if the $x_j$'s are B--splines. We assume the 
points $t_j$ are uniformly placed on the interval $[0,1]$ and define for 
$j=0,\dots,N-1$ the cubic B--splines

\begin{eqnarray*} 
S_j(t) = \frac{1}{4h^3}
\left\{ \begin{array}{ll}
(t-t_{j-2})^3                                    &,\ t\in [t_{j-2},t_{j-1}] \\
h^3+3h^2(t-t_{j-1})+3h(t-t_{j-1})^2-3(t-t_{j-1})^3 &,\ t\in [t_{j-1},t_j]   \\
h^3+3h^2(t_{j+1}-t)+3h(t_{j+1}-t)^2-3(t_{j+1}-t)^3 &,\ t\in [t_j,t_{j+1}]   \\
(t_{j+2}-t)^3                                      &,\ t\in [t_{j+1},t_{j+2}]
\end{array} \right.
\end{eqnarray*}

        In the formulation of our strategy we assumed the space $Y_h$ satisfy 
the condition in (\ref{Yh_wahl}). In order to rescue our choice of $Y_h$, we 
add to the system (\ref{nicht_lin_moment_probl_syst}) the following linear 
restriction to $\vphi$
\begin{equation}
  < \ \varphi , \ P_h A_1(x^0) - d P_h A_1(x^0)x^0 >_Y \ = \ 0.
\label{bg_acocha} \end{equation}
\noindent  Joining the equations in (\ref{nicht_lin_moment_probl_syst}) and 
(\ref{bg_acocha}), we have an overdetermined system with $N+1$ equations to 
solve, in order to determine the $N$ coefficients of $\vphi$. We can 
observe, that the matrix coefficients $a_{i,j}$ of this system vanish for 
$i > j+2$ and $i \not= N+1$.

\subsection{A linear elliptic Cauchy--problem}

        We begin with the definition of the linear operator $A: H^{\me}
(\Gamma_r) \mapsto H^{\me}(\Gamma_l)$

\vskip-.3cm \unitlength1cm
\centerline{\hfill
\begin{picture}(15,5)
\put(2.1,2.5){ \mbox{ $\left\{ \begin{array}{rl}
                \Delta w = 0,       & in\ \Omega \\
                w \      = \varphi, & at\ \Gamma_r\\
                w_\nu    = 0,       & at\ \Gamma_l\\
                w_\nu    = 0,       & at\ \Gamma_i
        \end{array} \right.$ } }
\put(11,2.5){\circle{2.0}}
\bezier64(9.5,2.5)(9.5,3.12)(9.94,3.56)    \bezier64(9.94,3.56)(10.38,4)(11,4)
\bezier64(12.5,2.5)(12.5,3.12)(12.06,3.56) \bezier64(12.06,3.56)(11.62,4)(11,4)
\bezier64(9.5,2.5)(9.5,1.88)(9.94,1.44)    \bezier64(9.94,1.44)(10.38,1)(11,1)
\bezier64(12.5,2.5)(12.5,1.88)(12.06,1.44) \bezier64(12.06,1.44)(11.62,1)(11,1)
\put(11,3.8){\line(0,1){0.4}} \put(11,1.2){\line(0,-1){0.4}}
\put(10.9,3.4){${\Omega}$}
\put(9.6,3.9){${\Gamma_l}$} \put(12.1,3.9){${\Gamma_r}$}
\put(7.4,3.0){$A(\varphi):=w_{|_{\Gamma_l}}$}
\put(10.4,1.4){${\Delta w=0}$}
\put(10.6,2.2){${}_{w_\nu=0}$}
\put(8.5,1.1){$w_\nu=0$}
\put(12.3,1.1){$w=\varphi$}
\end{picture}
\hfill }
\vskip-.6cm

\noindent  were $H^s$ are the Soblev spaces%
\footnote{For details see [Ad] or [DaLi].}
of index $s \in \R$ and $w$ is the $H^1(\Omega;\Delta)$--solution of the 
mixed boundary value problem on the left hand side.

        Note that solving the equation $A \, \vphi = f$ is equivalent of 
finding the trace $\vphi = w_{|_{\Gamma_r}}$ of the 
$H^1(\Omega;\Delta)$--solution of following elliptic Cauchy--problem:
\begin{eqnarray*}
        \left\{ \begin{array}{rl}
                \Delta w = 0,  & in\ \Omega \\
                w \      = f,  & at\ \Gamma_l\\
                w_\nu    = 0,  & at\ \Gamma_l\\
                w_\nu    = 0,  & at\ \Gamma_i
        \end{array} \right.
\end{eqnarray*}
        Given a Distribution $\mu \in H^{-\me}(\Gamma_r)$ we will use the 
Backus--Gilbert strategy to approximate the value $< \mu, \varphi>$ by 
$<\psi , f>$, were $\psi$ is the solution of
\begin{equation}
                               A^* \psi \ = \mu.
\label{adjun_gleichung} \end{equation}
        We can see, using integration by parts, that the adjoint operator 
of the restriction of $A$ to $H^{\me}_{00}(\Gamma_r)$ is the operator 
$A^\sharp: H^{\me}_{00}(\Gamma_l)' \mapsto H^{\me}_{00}(\Gamma_r)'$ defined by

\vskip-.4cm \unitlength1cm
\centerline{\hfill
\begin{picture}(15,5)
\put(2.1,2.5){ \mbox{ $\left\{ \begin{array}{rl}
                \Delta v = 0,       & in\ \Omega \\
                v \      = 0,       & at\ \Gamma_r\\
                v_\nu    = \psi,    & at\ \Gamma_l\\
                v_\nu    = 0,       & at\ \Gamma_i
        \end{array} \right.$ } }
\put(11,2.5){\circle{2.0}}
\bezier64(9.5,2.5)(9.5,3.12)(9.94,3.56)    \bezier64(9.94,3.56)(10.38,4)(11,4)
\bezier64(12.5,2.5)(12.5,3.12)(12.06,3.56) \bezier64(12.06,3.56)(11.62,4)(11,4)
\bezier64(9.5,2.5)(9.5,1.88)(9.94,1.44)    \bezier64(9.94,1.44)(10.38,1)(11,1)
\bezier64(12.5,2.5)(12.5,1.88)(12.06,1.44) \bezier64(12.06,1.44)(11.62,1)(11,1)
\put(11,3.8){\line(0,1){0.4}} \put(11,1.2){\line(0,-1){0.4}}
\put(10.9,3.4){${\Omega}$}
\put(9.6,3.9){${\Gamma_l}$} \put(12.1,3.9){${\Gamma_r}$}
\put(11.7,0.8){$A^\sharp(\psi):=v_{\nu|_{\Gamma_r}}$}
\put(10.4,1.4){${\Delta v=0}$}
\put(10.6,2.2){${}_{v_\nu=0}$}
\put(8.5,3.1){$v_\nu=\psi$}
\put(12.5,3.1)      {$v=0$}
\end{picture}
\hfill }
\vskip-.5cm

        For $\psi \in H^{\me}_{00}(\Gamma_l)'$, if $\vphi \in H^{\me} 
(\Gamma_r) \backslash H^{\me}_{00}(\Gamma_r)$, it's not true that
\begin{eqnarray*}
  \int_{\Gamma_l}{A(\varphi)\, \psi\ d\Gamma} \ = \
                      -\int_{\Gamma_r}{\varphi\, A^\sharp(\psi)\ d\Gamma}.
\end{eqnarray*}
\noindent  To correct this problem we need the following theorem.

\begin{theor} \label{korrektur_formel}

        For $a, b \in \R$ let $\eta_{a,b} \in C^\infty(\Gamma_r)$ be 
a function with $\eta_{a,b}(P_1) = a$ and $\eta_{a,b}(P_2) = b$, were 
$P_1$ and $P_2$ are the contact points between $\Gamma_r$ and $\Gamma_l$. If 
$V_{a,b}$ is the subspace of $H^{\me}(\Gamma_r)$ defined by
\begin{eqnarray*}
  V_{a,b} \ := \ \{ \vphi \in H^{\me}(\Gamma_r)\ /\ \eta_{a,b} - \vphi
                                            \in H^{\me}_{00}(\Gamma_r) \},
\end{eqnarray*}
\noindent  then for $\vphi_1, \vphi_2 \in V_{a,b}$ we have
\begin{eqnarray*}
             \int_{\Gamma_l}{A \vphi_1 \, \psi \ d\Gamma} \ + \
             \int_{\Gamma_r}{\vphi_1 \, A^\sharp \psi \ d\Gamma} \ = \
             \int_{\Gamma_l}{A \vphi_2 \, \psi \ d\Gamma} \ + \
             \int_{\Gamma_r}{\vphi_2 \, A^\sharp \psi \ d\Gamma}
\end{eqnarray*}
\noindent   for every $\psi$ in $H^{-\me}(\Gamma_l)$. \hfill 
\end{theor}

        A complete proof of this Theorem can be found in [Le]. A direct 
consequence of Theorem \ref{korrektur_formel} is that for $a, b \in \R$ one 
can defines over $H^{-\me}(\Gamma_l)$ the linear functional

\begin{eqnarray*}
  r_{a,b}(\psi) \ := \ < A \eta_{a,b} , \psi  > \ + \
                                             < \eta_{a,b} , A^\sharp \psi>
\end{eqnarray*}
\noindent  and obtain
\begin{eqnarray*}
  <A \vphi , \psi> \ = \ -< \vphi , A^\sharp \psi > \ + \ r_{a,b}(\psi),
\end{eqnarray*}

\noindent  for every $\vphi \in V_{a,b}$ and $\psi \in H^{-\me}(\Gamma_l)$.

        If we are able to find a $\psi \in H^{-\me}(\Gamma_l)$ that solves 
the equation
\begin{eqnarray*}
                          - A^\sharp \psi \ = \mu,
\end{eqnarray*}
\noindent  we can solve our reconstruction problem as before, using
\begin{eqnarray*}
        < \mu , \vphi > & = & - < A^\sharp \psi , \vphi> \\
                        & = & < \psi , A \vphi > \ - \ r_{a,b}(\psi) \\
                        & = & < \psi , f > \ - \ r_{a,b}(\psi).
\end{eqnarray*}
        We should observe that, if $\vphi \in H^{\me}_{00}(\Gamma_r)$, than 
$a = b = 0$ and $r_{a,b} \equiv 0$. In this case we have
\begin{eqnarray*}
                     < \mu , \vphi > \ = \ < \psi , f >.
\end{eqnarray*}

\section{Numerical results}

\subsection{The moment-problem}

        In this section we study the operator $A: L^2(0,1) \rightarrow 
L^2(0,1)$ defined in (\ref{nicht-lin_operator}) for $x^0(t) = t$. Let us start 
with the linear case, i.e., taking $\nu = 0$ in (\ref{nicht-lin_operator}).

        We generate different right hand sides by solving the 
direct problem $y = Ax$ for three functions
\begin{eqnarray*}
   x_a(t) = \left\{ \begin{array}{cl} t/2 & ,\ t \leq 1/2 \\
   t-\frac{1}{4} & ,\ t \geq 1/2 \end{array} \right.
\hskip 0.1cm , \hskip 0.6cm
   x_b(t) = \left\{\begin{array}{cl}  2t & ,\ t \leq 1/2 \\
   2-2t & ,\ t \geq 1/2 \end{array}\right.
\end{eqnarray*}
\noindent   and
\begin{eqnarray*}
   x_c(t) = \left\{\begin{array}{cl} 1 & ,\ \frac{1}{4}\leq t\leq
   \frac{3}{4} \\  0 & ,\ {\rm otherwise} \end{array}\right. .
\end{eqnarray*}

\noindent  Our grid is defined by $t_j = j/N,\ 0\leq j \leq N$. For the space 
$X_h$ we choose the B--spline basis corresponding to this grid. Our objective 
is to reconstruct the values of the different solutions $x_a$, $x_b$ and 
$x_c$ at the grid points $t_j$ and at the points $(t_{j+1} + t_j) /2,\ 
0 \leq j \leq N-1$. In Figure~1 we give the results for $N = 25$ and 
$N = 50$, when the exact right hand side $y$ is used.

\begin{figure}
\epsfysize9.2cm \epsfxsize14cm
\centerline{ \epsfbox{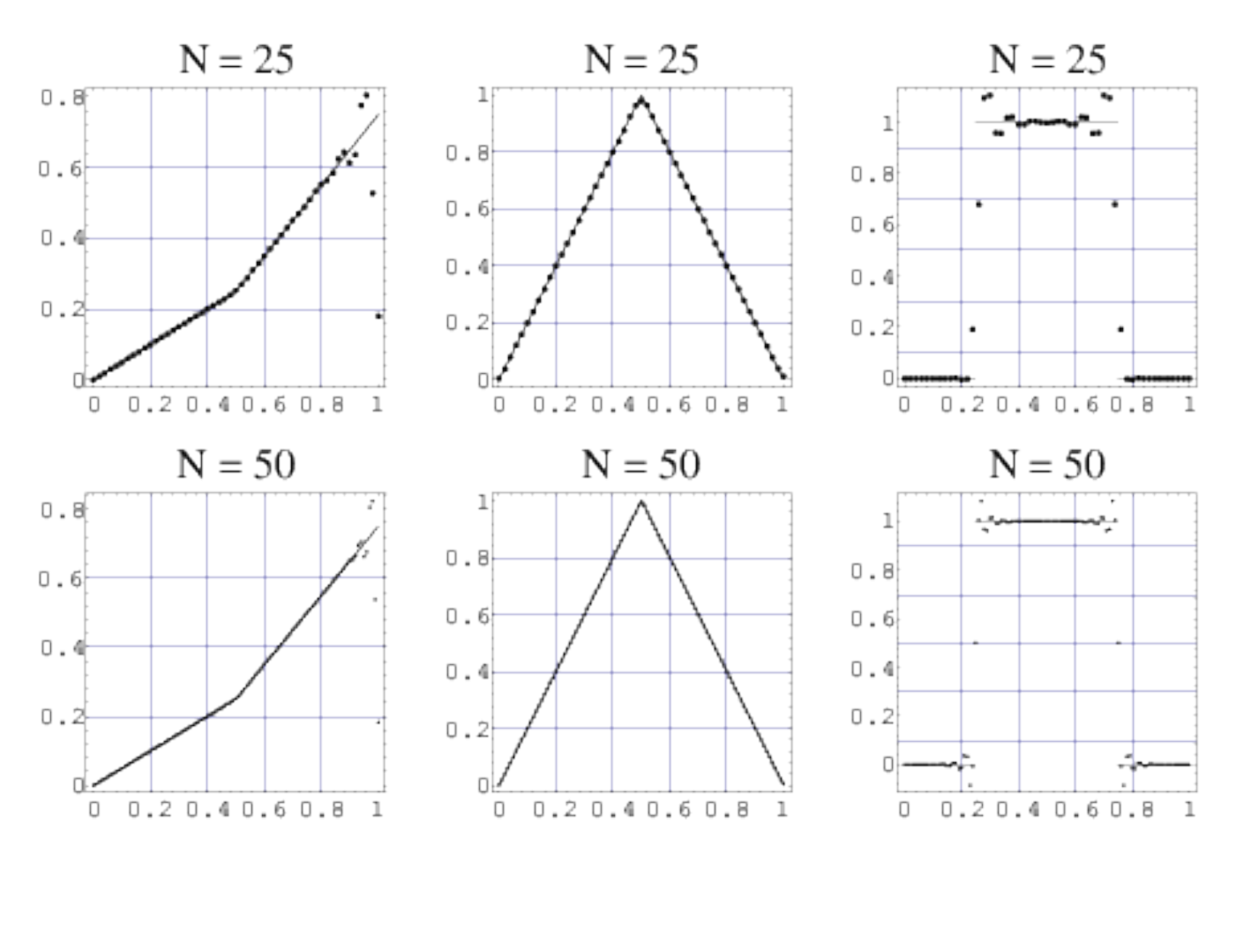} }
\vskip -1cm
\caption{}                                                           
\end{figure}

        In Figure~2, we show the reconstruction results for the linear 
operator and perturbed data. The system is solved for a right hand side 
$y_\eps$ generated by adding a $1\%$ random noise to the original $y$, 
i.e., $|y_j - y_{\eps,j}| \leq \frac{1}{100} y_j$.

\begin{figure}
\epsfysize3.8cm \epsfxsize14cm
\centerline{ \epsfbox{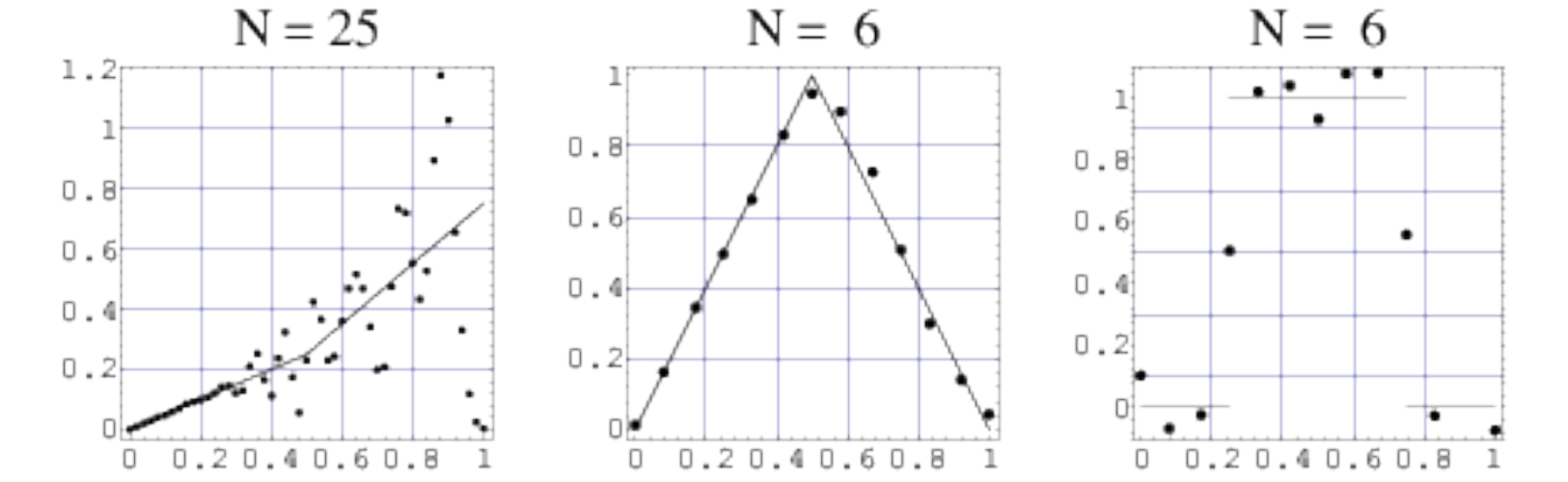} }
\caption{}                                                           
\end{figure}

        Next we analyze the reconstruction error at the point $t = 
\frac{1}{2}$ for exact data and the functions
\begin{eqnarray*}
   x_a(t) \ = \ 2t \hskip 0.6cm {\rm and} \hskip 0.6cm
   x_b(t) = \left\{\begin{array}{cl}  2t & ,\ t \leq 1/2 \\
   2-2t & ,\ t \geq 1/2 \end{array}\right. .
\end{eqnarray*}

        Analyzing Figure~3 we observe that the reconstruction is somehow 
better for even values of $N$. This can be explained by the existence of 
a B--spline centered at the point $t = \frac{1}{2}$ in the $X_h$--basis. A 
consequence of this is that the functional $\delta(\cdot - \frac{1}{2})$ will 
be better approximated in $X_h$ if $N$ is even.

\begin{figure}
\epsfysize4.4cm \epsfxsize12.2cm
\centerline{ \epsfbox{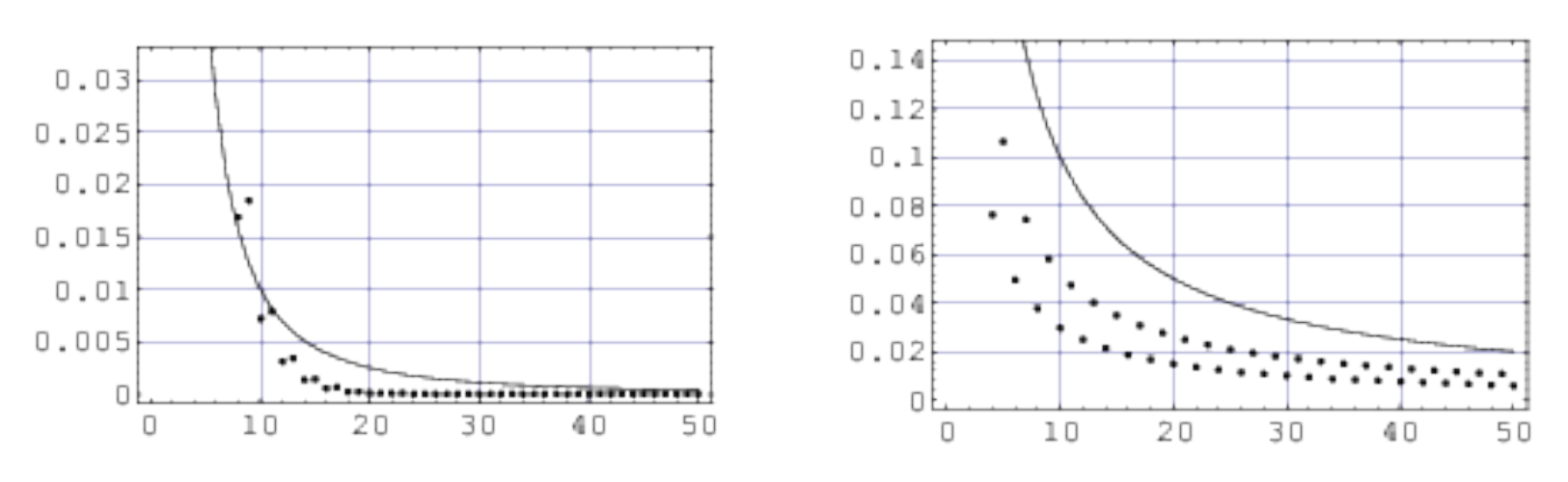} }
\vskip-5cm
\noindent \unitlength1cm \centerline{ \begin{picture}(14.4,5)
\put(3.5,3.2){${\frac{1}{x^2}}$} \put(9.8,3.4){${\frac{1}{x}}$}
\put(2.2,4.3){\small error for $x_a$} \put(8.4,4.3){\small error for $x_b$}
\put(7,0.8){N} \put(13.4,0.8){N}
\end{picture} }
\caption{}                                                           
\end{figure}

        Next we will analyze the operator $A$ in (\ref{nicht-lin_operator}) 
for small values of $\nu$. We use the same grid as before with $N = 25$ and 
try to reconstruct the polynomial $x^2$ at the points $t_j$ and $(t_{j+1} + 
t_j)/2$ using exact data. The results are shown in Figure~4.

\begin{figure}
\epsfysize4cm \epsfxsize14.2cm
\centerline{ \epsfbox{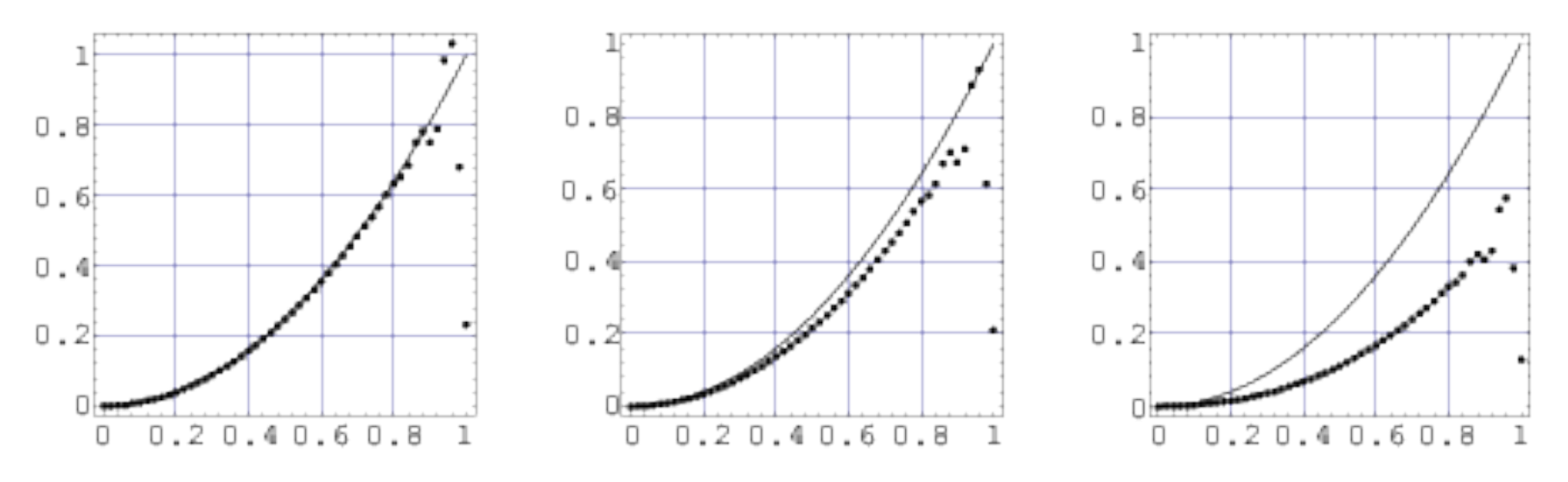} }
\vskip-4cm
\noindent \unitlength1cm \centerline{ \begin{picture}(14.2,4)
\put(2.0,4.1){$\nu=0.01$} \put(6.9,4.1){$\nu=0.1$} \put(11.8,4.1){$\nu=1$}
\end{picture} }
\caption{}                                                           
\end{figure}

        The next example in Figure~5 shows a reconstruction for $\nu = 0.01$ 
and exact data of the functions
\begin{eqnarray*}
   x_b(t) = \left\{\begin{array}{cl}  2t & ,\ t \leq 1/2 \\
   2-2t & ,\ t \geq 1/2 \end{array}\right.
\hskip 0.7cm {\rm and} \hskip 0.7cm
   x_c(t) = \left\{\begin{array}{cl} 1 & ,\ \frac{1}{4}\leq t\leq
   \frac{3}{4} \\  0 & ,\ {\rm otherwise} \end{array}\right. .
\end{eqnarray*}

\begin{figure}
\epsfysize3.9cm \epsfxsize13.8cm
\centerline{ \epsfbox{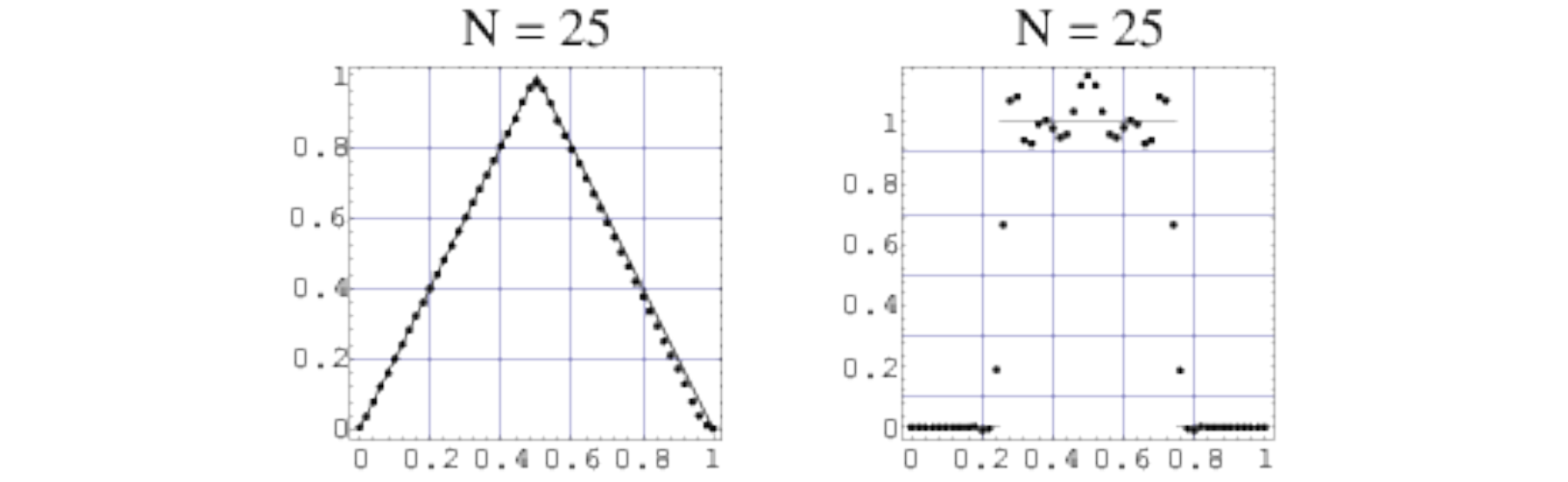} }
\bigskip
\caption{}                                                           
\end{figure}

\subsection{The elliptic Cauchy--problem}

        We will analyze the elliptic Cauchy--problem in an anulus $\Omega$ 
with inner radius $\frac{1}{2}$ and outer radius 1. Let us take the linear 
operator $A$ defined in \S 3.2. The problem we want to solve is, given a $\mu 
\in H^{-\me}(\Gamma_r)$, reconstruct the value of $< \mu , \vphi >$, where 
$\vphi \in H^{\me}(\Gamma_r)$ is the solution of the equation $A \vphi = f$. 
In order to generate consistent data $f$, we solve the direct problems for 
$\vphi_1(t) = (t - \frac{\pi}{2})^2$ and $\vphi_2(t) = \pi - 2|t - 
\frac{\pi}{2}|$, were $t \in [0,\pi]$.

        The formulation of this elliptic Cauchy--problem in $\Omega$ involves 
an extra difficulty: we are not able to characterize the space $Rg(A^\sharp)$. 
As we do want to have an element $\mu \in Rg(A^\sharp)$, we solve first the 
direct problem $\mu = A^\sharp \psi$ for a $\psi \in H^{-\me}(\Gamma_l)$. For 
this propose we chose $\psi \equiv 1$, solve the mixed boundary value problem
\begin{eqnarray*}
        \left\{ \begin{array}{rl}
                \Delta v = 0,     & in\ \Omega \\
                v \      = 0,     & at\ \Gamma_r\\
                v_\nu    = \psi,  & at\ \Gamma_l\\
                v_\nu    = 0,     & at\ \Gamma_i
        \end{array} \right.
\end{eqnarray*}
\noindent  and set $\mu = v_{\nu|_{\Gamma_r}} \in H^{-\me}(\Gamma_r)$. 

        According to the Backus--Gilbert strategy discussed in \S 3.2, 
the first thing to do is to solve the equation $-A^\sharp \psi = \mu$. To 
approximate the solution $\psi$, we use the iterative method described in 
[MaKo] (this iterative method is also extensively discussed in [Le]). The 
approximations $\psi_k$ are shown in Figure~6, were $k$ represents the 
iteration index. The grid node 0 represents the point $(0,-1)$ and the grid 
node 32 the point $(0,1)$ of $\Gamma_r$.

\begin{figure}
\epsfysize4cm \epsfxsize14.2cm
\centerline{ \epsfbox{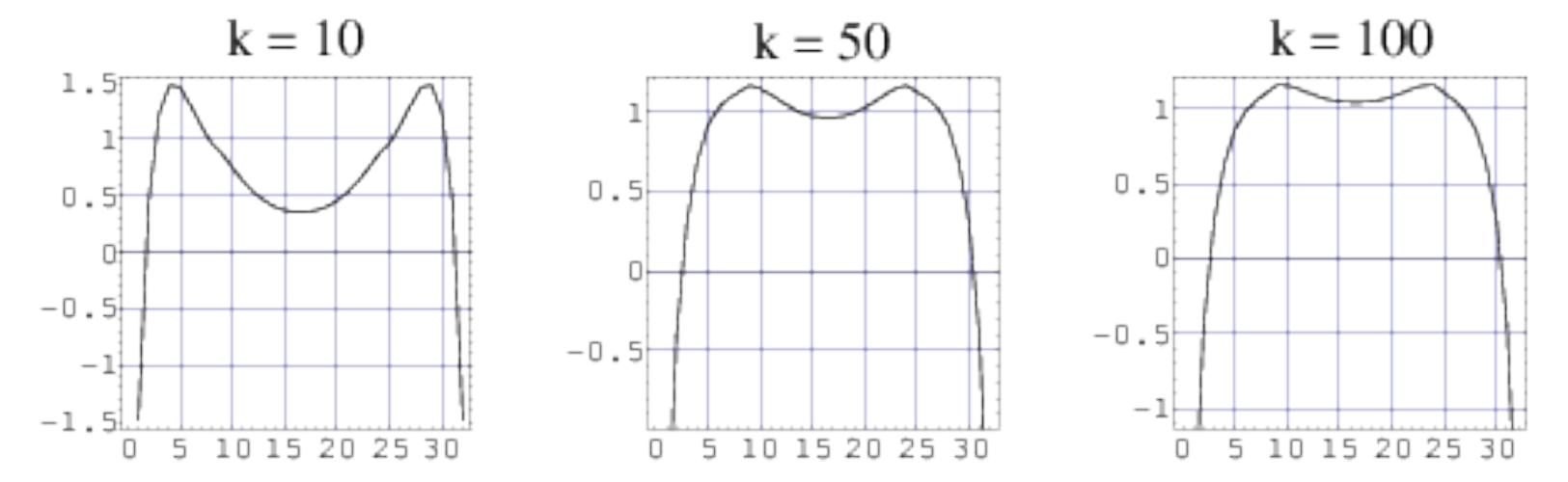} }
\vskip-4cm
\unitlength1cm \centerline{ \begin{picture}(14.2,4)
\put(4.6,0.3){$\Gamma_r$} \put(9.4,0.3){$\Gamma_r$} \put(14.1,0.3){$\Gamma_r$}
\end{picture} }
\caption{}                                                           
\end{figure}

        What we do next is to compare the values $< \mu , \vphi >$ with 
$< \psi , f > - \ r_{a,b}(\psi)$. The results are shown in Table~1 (note that 
$r_{a,b} \equiv 0$ for $\vphi = \vphi_2$).

{\rm
\begin{table} \caption{\mbox{ }\quad\quad\quad\quad\quad\quad\quad\quad\quad\quad
                               \quad\quad\quad\quad\quad\quad\quad\quad\quad\quad
                               \quad\quad\quad\quad\quad\quad\quad\quad\quad\quad
                               \quad\quad\quad\quad\quad\quad\quad\quad\quad\quad\mbox{}}
\begin{tabular}{@{}ccccc} 
\hline
{}  &  $< \mu , \vphi >$  &  $< \psi, f >$  &
$< \psi, f > - \ r_{a,b}(\psi)$ & relative error \\
\hline
$\vphi=\vphi_1$ &  $3.99577$ &  $5.84496$ &  $3.89191$ &  $0.01776$  \\
$\vphi=\vphi_2$ &  $0.52309$ &  $0.49903$ &  $0.49903$ &  $0.04599$  \\
\hline
\end{tabular}
\end{table} }

        Other numerical tests related to this specific Cauchy problem and to 
the validation of Theorem~\ref{korrektur_formel} can be found in [Le].

\section{Final remarks and Conclusions}

\mbox{} \indent
        $\bullet$ \ The numerical experiments show that one can obtain good 
approximations for $\mu = \delta$ in $Rg(dA^* P_h^*)$ if $A$ is the integral 
operator defined in \S 3.1 and the non-linearity in $A$ is small. In 
the non-linear case we can always improve an approximation defining a 
new $\tilde{x}^0$ as the B--spline interpolation of the evaluated values 
$x(t_j)$ and solving the new system
\begin{eqnarray*}
    < dA^*(\tilde{x}^0) \, P_h^* \vphi, \, x_j >_X \ = \ < \mu , x_j >_X.
\end{eqnarray*}
        Comparable and related results can be found in [Ch], [Hu], [Ki], 
[Lo], [LM1,2], [ScBe] and [Sn].

        $\bullet$ \ We should observe an unwanted Gibb's phenomenon in 
Figures~1a, 2a and 4. An explanation for this fact is that $\delta(\cdot) 
\in H^{-\me-\eps}$ for $\eps>0$ but $H^s_0([0,1]) \not\subseteq H^s([0,1])$ 
for $s > 1/2$. Thus the inner product $<\delta, x>_{L^2}$ will be in duality 
only if the boundary conditions $x(0)=x(1)=0$ are satisfied. \\
The same phenomenon can also be observed in Figure~5b, were the lack of 
regularity of the solution $x_c$ is now responsible for the effect.

        $\bullet$ \ If the operator $A$ is defined by the elliptic 
Cauchy problem in \S 3.2, we do not know, for an arbitrary set 
$\Omega$,  how to characterize the space $Rg(A^\sharp)$. But if some argument 
guarantees that the $\mu_i$'s are in $Rg(A^\sharp)$, we can proceed as in 
\S 4.2 and solve the Cauchy problems $A^\sharp \psi_i = \mu_i$ once 
for each $\mu_i$, in order to obtain the {\em observations}
\begin{eqnarray*}
                   < \mu_i , \vphi > \ = \ < \psi_i , f >
\end{eqnarray*}
\noindent  of $\vphi$, every time we have a different set of data $f$. Such 
$\mu_i$'s are also known in the literature as {\em sentinels} (see [Ch]).

        $\bullet$ \ When we analyzed the Cauchy--problem, we tried first to 
evaluate the reconstruction with $\mu = \delta$ and $\mu$ a 
$C^\infty$--mollifier. Using classical arguments (see [GiTr]) one can prove 
that no analytical solution exists in such cases when $\Omega$ has an 
analytical boundary. Our numerical results showed, that in this cases the 
equation $A^\sharp \psi = \mu$ has no solutions.

        $\bullet$ \ It is important to point out here the ill-posed nature of 
the involved reconstruction problems. Fredholm operators of the first kind 
are typically ill-posed [Gro]. What concern the elliptic Cauchy--problems, 
Hadamard elaborated an example with Cauchy data that converge uniformly to 
zero but the respective solutions become unbounded. The example follows:
\begin{eqnarray*}
\left\{  \begin{array}{rll}
  \Delta u_k = &\!\!\!0          &\!\!\!,\ (x,y)\in \Omega=(0,1)\times(0,1) \\
  u_k(x,0) =   &\!\!\!0          &\!\!\!,\ x \in (0,1) \\
  \frac{\partial}{\partial y} u_k(x,0) = &\!\!\!\vphi_k(x)&\!\!\!,\ x\in (0,1)
\end{array} \right.
\end{eqnarray*}
\noindent  were $\varphi_k = (\pi k)^{-1}sin(\pi k x)$. The respective 
solutions are
\begin{eqnarray*}
            u_k(x,y)\ =\ (\pi k)^{-2} sinh(\pi k y)\ sin(\pi k x).
\end{eqnarray*}
\indent
        $\bullet$ \ Our numerical experiments were realized on a IBM RISC 
6000/250 Work Station. It took some seconds to generate and solve the 
systems in \S 4.1 for $N = 50$. To evaluate the first 100 steps of the 
iterative method, in order to solve the Cauchy--Problem in \S 4.2, we needed 
about 30 minutes CPU-time (we used the finite element method on a grid with 
$\simeq$ 8000 nodes to solve each mixed BVP involved on the iterative method).

\end{document}